\documentclass[a4paper, 12pt, twoside]{article}
\usepackage{amsmath}
\usepackage{amsthm}
\usepackage{amsfonts}
\usepackage{pifont}
\usepackage{bbm}
\usepackage{mathrsfs}
\usepackage{fancyhdr}
\usepackage{graphicx}

\usepackage{psfrag}
\usepackage{time}
\usepackage{txfonts}

\usepackage{stmaryrd}
\usepackage{mathrsfs}
\usepackage{txfonts}
\usepackage{amsfonts}

\usepackage{indentfirst}

\newcommand{\xiaosan}{\fontsize{15pt}{22pt}\selectfont}
\newcommand{\sihao}{\fontsize{14pt}{21pt}\selectfont}

\newcommand{\xiaosi}{\fontsize{12pt}{18pt}\selectfont}

\usepackage{curves}

\numberwithin{equation}{section}

\newtheorem{theorem}{ {Theorem}}[section]

\newtheorem{remark} {   {Remark}}[section]

\newtheorem{lemma} {  {Lemma}}[section]
\newtheorem{definition}{  {Definition}}[section]

\begin{document}
\setlength{\parindent}{2em}
\newpage
\fontsize{12}{22}\selectfont\thispagestyle{empty}
\renewcommand{\headrulewidth}{0pt}
 \lhead{}\chead{}\rhead{} \lfoot{}\cfoot{}\rfoot{}
\noindent

\title{\bf Anderson  localization  for  the completely resonant  phases }
\author{{Wencai  Liu and Xiaoping Yuan*}\\
{\em\small School of Mathematical Sciences}\\
{\em\small  Fudan University}\\
{\em\small  Shanghai 200433, People's Republic of China}\\
{\small 12110180063@fudan.edu.cn}\\
{\small *Corresponding author: xpyuan@fudan.edu.cn}}
\date{}
\maketitle

\renewcommand{\baselinestretch}{1.2}
\large\normalsize

\begin{abstract}
For the
almost Mathieu operator $ (H_{\lambda,\alpha,\theta}u) (n)=u(n+1)+u(n-1)+ \lambda v(\theta+n\alpha)u(n)$,
 Avila and  Jitomirskaya guess that for every phase  $ \theta \in \mathscr{R} \triangleq\{\theta\in \mathbb{R}\;|   \;   2\theta + \alpha \mathbb{Z} \in \mathbb{Z}\}$,
 $H_{\lambda,\alpha,\theta}$ satisfies  Anderson localization  if $ |\lambda| >  e^{ 2 \beta}$. In the present paper, we show that  for every phase  $ \theta \in \mathscr{R}  $,
 $H_{\lambda,\alpha,\theta}$ satisfies  Anderson localization   if $ |\lambda| >  e^{  7 \beta}$.

\end{abstract}

\setcounter{page}{1} \pagenumbering{arabic}\topskip -0.82in
\fancyhead[LE]{\footnotesize  Introduction}
\section{\xiaosan \textbf{Introduction}}
   The almost Mathieu operator (AMO) is the   quasi-periodic   Schr\"{o}dinger operator  on  $   \ell^2(\mathbb{Z})$:
 \begin{equation}\label{G11}
 (H_{\lambda,\alpha,\theta}u) (n)=u(n+1)+u(n-1)+ \lambda v(\theta+n\alpha)u(n),   \text{ with }  v(\theta)=2\cos2\pi \theta,
 \end{equation}
where $\lambda$ is the coupling, $\alpha $ is the frequency, and $\theta $ is the phase.
 \par
 AMO is the most studied quasi-periodic Schr\"{o}dinger operator, arising naturally as a physical model
  (see \cite{Las} for a recent historical account and  for the physics background).
  \par
  We say phase $\theta\in \mathbb{R}$ is completely resonant with respect to frequency $\alpha$, if $\theta\in \mathscr{R} \triangleq\{\theta\in \mathbb{R}\;|   \;   2\theta + \alpha \mathbb{Z} \in \mathbb{Z}\}$.
 \par
Anderson localization (i.e.,  only pure point spectrum with exponentially decaying eigenfunctions) is not only
meaningful in physics, but also relates to reducibility  for   Aubry dual model(see \cite{GJLS}). In particular, Anderson localization for completely resonant phases  is crucial to  describe      open gaps of
$  \Sigma_{\lambda,\alpha}$( the spectrum of $H_{\lambda,\alpha,\theta}$ is independent of $\theta$ for $\alpha\in \mathbb{R}\backslash\mathbb{Q}$, denoted by $  \Sigma_{\lambda,\alpha}$).  See   $\cite{AJ2}$, \cite{LIU1}   and \cite{P2} for details.
 \par
It is well known  that $H_{\lambda ,\alpha,\theta}$  has purely absolutely continuous spectrum for $\alpha \in \mathbb{Q}$ and all $\lambda$.
This implies $H_{\lambda ,\alpha,\theta}$  can not satisfy Anderson localization for all $\alpha \in \mathbb{Q}$. Thus
 we always assume  $\alpha\in \mathbb{R}\backslash \mathbb{Q}$  in the present paper.
\par
   The following notions are essential in the study of equation  (\ref{G11}).
  \par
 We say $\alpha \in \mathbb{R}\backslash \mathbb{Q} $ satisfies a Diophantine condition $\text{DC}(\kappa,\tau)$ with $\kappa>0$ and $\tau>0$,
if
$$ ||k\alpha||_{\mathbb{R}/\mathbb{Z}}>\kappa |k|^{-\tau}  \text{ for   any } k\in \mathbb{Z}\setminus \{0\},$$
where $||x||_{\mathbb{R}/\mathbb{Z}}=\min_{\ell \in \mathbb{Z}}|x-\ell| $.
Let $\text{DC}=\cup_{\kappa>0,\tau>0}\text{DC}( \kappa,\tau)$. We say $\alpha$  satisfies   Diophantine condition, if $\alpha\in \text{DC}$.
 \par
Let
 \begin{equation}\label{G12}
  \beta= \beta(\alpha)=\limsup_{n\rightarrow\infty}\frac{\ln q_{n+1}}{q_n},
 \end{equation}
 where $ \frac{p_n}{q_n} $ is  the continued fraction approximants   to $\alpha$. Notice that $\beta(\alpha)=0$ for $\alpha\in DC$.
 \par

 Avila and Jitomirskaya guess that for any   completely resonant phase $\theta$, $H_{\lambda,\alpha,\theta}$ satisfies Anderson localization  if $ |\lambda| >  e^{  2\beta}$ (Remark 9.1,  $\cite{AJ1}$).
 Jitomirskaya-Koslover-Schulteis
    proves this  for  $\alpha\in DC  $   \cite{JKS}, more concretely, for $\alpha\in DC$, $H_{\lambda,\alpha,\theta}$ satisfies Anderson localization  if   $\theta\in \mathscr{R}$ and $ |\lambda| > 1$.
   In $ \cite{AJ2}$,  Avila and Jitomirskaya
   firstly develop    a quantitative version of  Aubry duality.
   By the way, they  obtain that for $\alpha$ with $\beta(\alpha)=0$,
    $H_{\lambda,\alpha,\theta}$ satisfies Anderson localization  if $\theta\in \mathscr{R}$  and  $ |\lambda| > 1$.
     The present authors extend the quantitative version of  Aubry duality to all $\alpha$ with $\beta(\alpha)<\infty$,
     and show that  for all $\alpha$ with $\beta (\alpha)<\infty$,
    $H_{\lambda,\alpha,\theta}$ satisfies Anderson localization  if $\theta\in \mathscr{R}$  and  $ |\lambda| > e^{C\beta}$, where
    $C$ is a large absolute constant\cite{LIU1}.
     In the present  paper,  we   give a definite quantitative description  about the constant $C$,
and obtain  the following theorem.
\begin{theorem}\label{Main Theorem}
 For $\alpha\in  \mathbb{R}\backslash \mathbb{Q}$ with $\beta (\alpha)<\infty$, the  almost Mathieu operator
    $H_{\lambda,\alpha,\theta}$ satisfies Anderson localization  if $\theta\in \mathscr{R}$  and  $ |\lambda| > e^{7\beta}$, where
     $ \mathscr{R} \triangleq\{\theta\in \mathbb{R}\;|   \;   2\theta + \alpha \mathbb{Z} \in \mathbb{Z}\}$.
\end{theorem}
\begin{remark}
 Avila-Jitomirskaya thinks  that $H_{\lambda,\alpha,0}$  does not display  Anderson localization  if $   |\lambda| \leq  e^{  2\beta}$ (Remark 5.2,  $\cite{AJ1}$), which  is still open.  Clearly,     $0\in \mathscr{R}$.
\end{remark}
 Avila and  Jitomirskaya guess that for a.e. $\theta$, $H_{\lambda,\alpha,\theta}$ satisfies Anderson localization  if $ |\lambda| >  e^{  \beta} $(Remark 9.2,  $\cite{AJ1}$), and they   establish  this for $ |\lambda| >  e^{\frac{16}{9} \beta}$. This  result
 has been extended to     regime $ |\lambda| >  e^{\frac{3}{2} \beta}$ by the present authors\cite{LIU}. More precisely, there exists
 a Lebesgue zero-measure $B$ such that
 $H_{\lambda,\alpha,\theta}$ satisfies Anderson localization  if $\theta\notin B$ and $ |\lambda| >  e^{\frac{3}{2}  \beta}$.
 Unfortunately, $\mathscr{R}\subset B $.   In the present paper,  we make   some adjustment such that  the discussion in \cite{AJ1}
 and  \cite{LIU} can be applied to completely resonant phase $\theta$.
 \par
 The present paper is organized as follows:
 \par
In \S 2,  we give some preliminary notions and facts which are taken from  Avila-Jitomirskaya  \cite{AJ1} or Bourgain \cite{B2}. In \S3, we set up the regularity of  non-resonant $y$.
 In \S4, we set up the regularity of  resonant $y$.
 In \S5, we give the proof of Theorem \ref{Main Theorem} by the regularity of $y$ and  block resolvent expansion.

   \section{\xiaosan \textbf{ Preliminaries}}

  It is well known that  Anderson localization for a self-adjoint operator  $H$ on $\ell^2(\mathbb{Z})$ is equivalent to  the following statements.
  \par
  Assume $\phi$ is an extended state of $H$, i.e.,
  \begin{equation}\label{G21}
    H\phi=E\phi  \text{ with  } E\in \Sigma(H) \text{ and }  |\phi(k)|\leq (1+|k|)^C,
  \end{equation}
  where $\Sigma(H)$ is the spectrum of   self-adjoint   operator $H$.
  Then there exists some constant $c>0$ such that
  \begin{equation}\label{G22}
    | \phi(k)|< e^{-c|k|} \text{ for } k\rightarrow \infty.
  \end{equation}
 \par
 The above statements can be proved by Gelfand-Maurin Theorem.
 See \cite{BER} for the proof of    continuous-time Schr\"{o}dinger operator.
The proof of  discrete     Schr\"{o}dinger operator is similar, see   \cite{LIU2} for example.

If   $ \alpha $ satisfies  $ \beta(\alpha) =0$,   Theorem   \ref{Main Theorem}  has been proved by  Avila-Jitomirskaya  $\cite{AJ2} $,
 which we have mentioned in \S 1. Thus in the present  paper, we fix $ \alpha\in \mathbb{R}\backslash\mathbb{Q}$
  such that $ 0<\beta(\alpha) <\infty$.  Unless  stated otherwise,   we always  assume $\lambda>e^{7\beta}$ ( for $\lambda<-e^{7\beta}$,  notice that   $H_{\lambda,\alpha,\theta}=H_{-\lambda,\alpha,\theta+\frac{1}{2}}$), and
 $E\in \Sigma_{\lambda,\alpha}$. Since this does not   change  any of the statements,  sometimes the dependence of parameters $E,\lambda,\alpha, \theta$ will be ignored in the following.
  \par
  Given      an extended state $\phi$ of $H_{\lambda,\alpha,\theta}$, without loss of generality assume $\phi(0)=1$.
  Our objective is to show that there exists some $c>0$ such that
  $$| \phi(k)|< e^{-c|k|} \text{ for } k\rightarrow \infty.$$

\par
Let us denote
$$ P_k(\theta)=\det(R_{[0,k-1]}(H_{\lambda,\alpha,\theta}-E) R_{[0,k-1]}).$$
Following $ {\cite{JKS}}$, $P_k(\theta)$ is an even function of $ \theta+\frac{1}{2}(k-1)\alpha$  and can be written as a polynomial
      of degree $k$ in $\cos2\pi (\theta+\frac{1}{2}(k-1)\alpha )$ :
     \begin{equation}\label{G23}
      P_k(\theta)=\sum _{j=0}^{k}c_j\cos^j2\pi (\theta+\frac{1}{2}(k-1)\alpha)    \triangleq  Q_k(\cos2\pi  (\theta+\frac{1}{2}(k-1)\alpha)).
     \end{equation}

      Let $A_{k,r}=\{\theta\in\mathbb{R} \;|\;Q_k(\cos2\pi   \theta   )|\leq e^{(k+1)r}\} $ with $k\in \mathbb{N}$ and $r>0$.
\begin{lemma}$ (\text{p}. 16, \cite{AJ1})$ \label{Le21}
The following inequality  holds
    \begin{equation}\label{G24}
     \lim_{k\rightarrow\infty}\sup_{\theta\in \mathbb{R}} \frac{1}{k} \ln | P_k(\theta)|\leq \ln \lambda.
    \end{equation}
\end{lemma}

    \par
    By Cramer's rule (p. 15, $\cite{B2}$)  for given  $x_1$ and $x_2=x_1+k-1$, with
     $ y\in I=[x_1,x_2] \subset \mathbb{Z}$,  one has
     \begin{eqnarray}
       |G_I(x_1,y)| &=&  \left| \frac{P_{x_2-y}(\theta+(y+1)\alpha)}{P_{k}(\theta+x_1\alpha)}\right|,\label{G25}\\
       |G_I(y,x_2)| &=&\left|\frac{P_{y-x_1}(\theta+x_1\alpha)}{P_{k}(\theta+x_1\alpha)} \right|.\label{G26}
     \end{eqnarray}
By Lemma \ref{Le21}, the numerators in  (\ref{G25}) and (\ref{G26}) can be bounded uniformly with respect to $\theta$. Namely,
for any $\varepsilon>0$,
\begin{equation}\label{G27}
    | P_n(\theta)|\leq e^{(\ln \lambda+\varepsilon)n}
\end{equation}
for $n$ large enough.
\begin{definition}\label{Def21}
Fix $t > 0$. A point $y\in\mathbb{Z}$ will be called $(t,k)$-regular if there exists an
interval $[x_1,x_2]$  containing $y$, where $x_2=x_1+k-1$, such that
\begin{equation}\label{G28}
  | G_{[x_1,x_2]}(y,x_i)|<e^{-t|y-x_i|} \text{ and } |y-x_i|\geq \frac{1}{7} k , \text{ for }i=1,2;
\end{equation}
otherwise, $y$ will be called $(t,k)$-singular.
\end{definition}
It is  easy to check that (p. 61, $\cite{B2}$)
 \begin{equation}\label{G29}
   \phi(y)= -G_{[x_1 ,x_2]}(x_1,y ) \phi(x_1-1)-G_{[x_1 ,x_2]}(y,x_2) \phi(x_2+1),
 \end{equation}
 where  $ y\in I=[x_1,x_2] \subset \mathbb{Z}$.
Our strategy is to establish the $(t,  k(y) )$-regular of $y$,  then localized property is easy to
obtain by   $(\ref{G29})$ and the block resolvent expansion.
       \begin{definition}
     We  say that the set $\{\theta_1, \cdots ,\theta_{k+1}\}$ is $ \epsilon$-uniform if
      \begin{equation}\label{G210}
        \max_{ x\in[-1,1]}\max_{i=1,\cdots,k+1}\prod_{ j=1 , j\neq i }^{k+1}\frac{|x-\cos2\pi\theta_j|}
        {|\cos2\pi\theta_i-\cos2\pi\theta_j|}<e^{k\epsilon}.
      \end{equation}
     \end{definition}
      \begin{lemma}\label{Le22}$(\text{Lemma 9.3 },\cite{AJ1})$
      Suppose  $\{\theta_1, \cdots ,\theta_{k+1}\}$ is  $ \epsilon_1$-uniform. Then there exists some $\theta_j$ in set  $\{\theta_1, \cdots ,\theta_{k+1}\}$ such that
     $\theta_j\notin A_{k,\ln\lambda-\epsilon}$ if    $ \epsilon>\epsilon_1$ and $ k$
      is sufficiently large.
      \end{lemma}
   Assume without loss of generality that $y>0$. Define $b_n=q_n^{8/9}$, where $q_n$ is given by (\ref{G12}), and  find $n$ such that $b_n\leq y<b_{n+1}$.
        We will distinguish two cases:
         \par
        (i)   $|y-\ell q_n|\leq b_n$ for some $\ell\geq1$,  called resonance.
          \par
        (ii)     $|y-\ell q_n|> b_n$ for all $\ell\geq0$, called  non-resonance.
\par
Next, we will  establish the  regularity  for resonant and non-resonant  $y$ respectively.
Given a phase  $ \theta\in \mathscr{R}$, there exists some $p\in \mathbb{Z}$ such that $ 2\theta-p\alpha\in \mathbb{Z} $.
Without loss of generality,    assume $p\leq 0$ below.

\section{ Regularity for non-resonant  $y$}
In this section, we will set up the regularity for non-resonant  $y$, for this reason, we give some lemmata  first.
Note that  $C$ is a large absolute constant  below,  which may change through the arguments, even when appear in the same formula.
For simplicity,  we    replace $I=[x_1,x_2]\cap \mathbb{Z}$  with  $I=[x_1,x_2]$.

   \begin{lemma}$ (\text{Lemma } 9.7, \cite{AJ1})$\label{Le31}
Let $\alpha\in \mathbb{R}\backslash \mathbb{Q}$, $x\in\mathbb{R}$ and $0\leq \ell_0 \leq q_n-1$ be such that
$ |\sin\pi(x+\ell_0\alpha)|=\inf_{0\leq\ell\leq q_n-1}    |\sin\pi(x+\ell \alpha)|$, then for some absolute constant $C > 0$,
\begin{equation}\label{G31}
    -Cq_n\leq \sum _{\ell=0,\ell\neq \ell_0}^{q_n-1} \ln|\sin\pi(x+\ell\alpha )|+(q_n-1)\ln2\leq  Cq_n,
\end{equation}
where
  $ q_n  $ is  given by (\ref{G12}).
  \par
Recall that     $\{q_n\}_{n\in \mathbb{N}}$ is the sequence of best denominators of irrational number $\alpha$,
since it satisfies
\begin{equation}\label{G32}
\forall 1\leq k <q_{n+1}, \| k\alpha\|_{\mathbb{R}/\mathbb{Z}}\geq ||q_n\alpha||_{\mathbb{R}/\mathbb{Z}},
\end{equation}
 Moreover, we also have the following estimate.
\begin{equation}\label{G33}
      \frac{1}{2q_{n+1}}\leq\Delta_n\triangleq \|q_n\alpha\|_{\mathbb{R}/\mathbb{Z}}\leq\frac{1}{q_{n+1}}.
\end{equation}
\end{lemma}
\par
 Now that  $y$ is non-resonant. Without loss of generality, let $y=mq_n+y_0$ with  $m\leq\frac{q_{n+1}^{ 8/9} }{q_n} $ and $q_n^{ 8/9}\leq y_0\leq \frac{q_n}{2}$.
 Let $s\in\mathbb{N}$   be
the largest positive integer  such that $4sq_{n-1}-p+1\leq    y_0$.
 Notice that $8sq_{n-1}<q_n$.
 \par
Set $I_1$, $I_2$ as follows,
$$I_1=[- 2sq_{n-1} , -1]$$ and
$$I_2=[mq_n+y_0- 2sq_{n-1}  ,mq_n+y_0+2sq_{n-1} -1].$$
The set $\{\theta_j\}_{j\in I_1\cup I_2}$
consists of $ 6sq_{n-1}$ elements, where $\theta_j=\theta+j\alpha $
    and $ j$ ranges through $I_1\cup I_2$.
\begin{lemma}\label{Le32}
 For any  $\varepsilon>0$,
the set $\{\theta_j\}_{j\in I_1\cup I_2}$ is $-2\ln (s /q_n)/q_{n-1}+\varepsilon$-uniform
  if $n$ is  sufficiently large.
\end{lemma}
\textbf{Proof:} We will first estimate numerator  in ($\ref{G210}$).  In ($\ref{G210}$), let $x=\cos2\pi a$ and take the logarithm, one has
   $$ \sum _{j \in I_1  \cup I_2,j\neq i}\ln|\cos2\pi a-\cos2\pi \theta_{j}| \;\;\;\;\;\;\;\;\;\;\;\;\;\;\;\;\;\;\;\;\;\;\;\;\;\;\;\;\;\;\;\;\;\;\;\;\;\;\;\;\;\;\;\;\;\;\;\;\;\;\;\;\;\;\;\;\;\;$$
$$\;\;\;\;\;\;\;\;\;\;\;\;\;\;\;=\sum_{j\in I_1  \cup I_2,j\neq i}\ln|\sin\pi(a+\theta_{j})|+\sum_{j \in I_1  \cup I_2,j\neq i}\ln |\sin\pi(a-\theta_{j})|
+(6sq_{n-1}-1)\ln2  $$
\begin{equation}\label{G34}
    =\Sigma_{+}+\Sigma_-+(6sq_{n-1}-1)\ln2,   \;\;\;\;\;\;\;\;\;\;\;\;\;\;\;\;\;\;\;\;\;\; \;\;\;\;\;\;\;\;\;\;\;\;\;\;\;\;\;\;\;\;\;\;\;\; \;
\end{equation}
 where
 \begin{equation}\label{G35}
   \Sigma_{+}=\sum_{j \in I_1  \cup I_2,j\neq i}\ln |\sin\pi(a+  \theta_{j}  )|,
 \end{equation}
 and
  \begin{equation}\label{G36}
     \Sigma_-=\sum_{j \in I_1  \cup I_2,j\neq i}\ln |\sin\pi ( a-\theta_{j})|.
 \end{equation}
Both $\Sigma_+$ and $\Sigma_-$ consist of $6s$ terms of the form of $(\ref{G31})$, plus 6s terms of the form
\begin{equation}\label{G37}
    \ln\min_{j=0,1,\cdots,q_n-1}|\sin\pi(x+j\alpha)|,
\end{equation}
minus $\ln|\sin\pi(a\pm\theta_i)|$. Since     there exists a interval of
length $q_n$ in sum of (\ref{G35}) ( or (\ref{G36}) ) containing $i$, thus the minimum  over this interval is not
 more than $\ln|\sin\pi(a\pm\theta_i)|$ (by the minimality).
       Thus,  using  $(\ref{G31})$ 6s times of $\Sigma_{+}$ and $\Sigma_{-}$ respectively,  one has
 \begin{equation}\label{G38}
   \sum_{j \in I_1  \cup I_2,j\neq i}\ln|\cos2\pi a-\cos2\pi \theta_{j}|\leq-6sq_{n-1}\ln2+Cs\ln q_{n-1}.
\end{equation}
The estimate of  the denominator of  ($\ref{G210}$) requires a bit more work.
 Without loss of generality, assume   $i\in I_1$.
 \par
 In $(\ref{G34})$, let $a=\theta_i$,
we obtain
   $$ \sum _{j \in I_1  \cup I_2,j\neq i}\ln|\cos2\pi \theta_i-\cos2\pi \theta_j| \;\;\;\;\;\;\;\;\;\;\;\;\;\;\;\;\;\;\;\;\;\;\;\;\;\;\;\;\;\;\;\;\;\;\;\;\;\;\;\;\;\;\;\;\;\;\;\;\;\;\;\;\;\;\;\;$$
$$\;\;\;\;\;\;\;\;\;\;\;\;\;\;\;\;\;\;\;\;\;=\sum_{j \in I_1  \cup I_2,j\neq i}\ln|\sin\pi(\theta_i+\theta_j)|+\sum_{j \in I_1  \cup I_2,j\neq i}\ln |\sin\pi(\theta_i-\theta_j)|
+(6sq_{n-1}-1)\ln2  $$
\begin{equation}\label{G39}
    =\Sigma_{+}+\Sigma_-+(6sq_{n-1}-1)\ln2,  \;\;\;\;\;\;\;\;\;\;\;\;\;\;\;\;\;\;\;\;\;\;\;\; \;\;\;\;\;\;\;\;\;\;\;\;\;\;\;\;\;\;
\end{equation}
 where
  \begin{equation}\label{G310}
   \Sigma_{+}=\sum_{j \in I_1  \cup I_2,j\neq i}\ln |\sin\pi(2\theta+ (i+j) \alpha)|.
 \end{equation}
 and
\begin{equation}\label{G311}
     \Sigma_-=\sum_{j \in I_1  \cup I_2,j\neq i}\ln |\sin\pi( i-j)\alpha|.
 \end{equation}
 We first estimate $\Sigma_+ $.
 Set
$J_1=[- 2s,-1]$ and
$J_2=[0 ,4s-1 ]$, which are two adjacent disjoint intervals of length
$2s$ and $4s$ respectively.
 Then $I_1\cup
I_2$ can be represented as a disjoint union of segments $B_j,\;j\in
J_1\cup J_2,$ each of length $q_{n-1}$.
Applying (\ref{G31}) on  each  $B_j$, we
obtain
\begin{equation}\label{G312}
\Sigma_+ > -6sq_{n-1}\ln 2+
\sum_{j\in J_1\cup J_2 }\ln
|\sin  \pi\hat \theta_j|-Cs\ln q_{n-1}-\ln|\sin2\pi (\theta+i\alpha)|,
\end{equation}
where
\begin{equation}\label{G313}
|\sin  \pi\hat \theta_j|=\min_{\ell \in B_j}|\sin  \pi
(2\theta +(\ell +  i)\alpha )|.
\end{equation}
We now start to estimating (\ref{G313}). Noting that $2\theta +(\ell +  i)\alpha\in (\ell +  i+p)\alpha +\mathbb{Z}$,
together with the construction of $I_1$ and $I_2$, one has
\begin{equation}\label{G314}
   2\theta +(\ell +  i)\alpha=mq_n\alpha+r_1\alpha\;\mod \mathbb{Z}
\end{equation}
or
\begin{equation}\label{G315}
   2\theta +(\ell +  i)\alpha=  r_2\alpha \mod \mathbb{Z},
\end{equation}
where $1\leq | r_i|<q_n$, $i=1,2$.
By (\ref{G32}) and  (\ref{G33}), we have
\begin{eqnarray}
\nonumber
  \min_{\ell \in I_1\cup I_2} ||2\theta +(\ell +  i)\alpha||_{\mathbb{R}/\mathbb{Z}}&\geq&  ||r_i\alpha||_{\mathbb{R}/\mathbb{Z}}-\frac{\Delta_{n-1}}{2} \\
  \nonumber
    &\geq&  \Delta_{n-1}-\frac{\Delta_{n-1}}{2} \\
    &\geq&  \frac{\Delta_{n-1}}{2},\label{G316}
\end{eqnarray}
since $ || mq_n\alpha||_{\mathbb{R}/\mathbb{Z}}\leq \frac{q_{n+1}^{ 8/9} }{q_n} \Delta_n\leq\frac{\Delta_{n-1}}{2}$.
\par
Next we will estimate $\sum_{j\in J_1  }\ln
|\sin  \pi\hat \theta_j|$.
Assume that $\hat \theta_{j+1}=\hat \theta_j+ q_{n-1}  \alpha$ for
every $j,j+1 \in J_1$.  Applying   the Stirling formula and (\ref{G316}), one has
 \begin{eqnarray}
  \nonumber
   \sum_{j\in J_1}\ln
|\sin 2\pi\hat\theta_j|&>&   2\sum_{j=1}^{s}\ln \frac
{j\Delta_{n-1}} {C}\\
    &>&  2s\ln\frac s{q_n}-C s.\label{G317}
 \end{eqnarray}
 \par
In the other case, decompose $J_1$ in maximal intervals $T_\kappa$
such that for $j,j+1 \in
T_\kappa$ we have $\hat \theta_{j+1}=\hat \theta_j+ q_{n-1}
\alpha$.  Notice that the boundary points of an interval $T_\kappa$
are either boundary points of $J_1$ or satisfy
$\| \hat \theta_j\|_{\mathbb{R}/\mathbb{Z}}+\Delta_{n-1} \geq
\frac {\Delta_{n-2}} {2}$.
This follows from the fact that if
$0<|z|<q_{n-1}$,  then $\|  \hat \theta_j+q_{n-1}\alpha\|_{\mathbb{R}/\mathbb{Z}}\leq \|  \hat \theta_j\|_{\mathbb{R}/\mathbb{Z}}+\Delta_{n-1}$,  and  $\|  \hat \theta_j+(z+q_{n-1})\alpha \|_{\mathbb{R}/\mathbb{Z}}\geq   \|z \alpha\|_{\mathbb{R}/ \mathbb{Z}}-\|  \hat \theta_j+q_{n-1}\alpha\|_{\mathbb{R}/\mathbb{Z}}
\geq \Delta_{n-2}-\|  \hat \theta_j\|_{\mathbb{R}/\mathbb{Z}}-\Delta_{n-1}$.
Assuming $T_\kappa \neq J_1$, then there exists $j \in
T_\kappa$ such that $\|  \hat \theta_j\|_{\mathbb{R}/\mathbb{Z}}\geq
\frac {\Delta_{n-2}} {2}- \Delta_{n-1} $.
\par
If $T_\kappa$  contains  some $j$ with $ \|  \hat
\theta_j\|_{\mathbb{R}/\mathbb{Z}}<\frac {\Delta_{n-2}} {10} $, then
 \begin{eqnarray}
 \nonumber
   |T_\kappa| &\geq & \frac{\frac {\Delta_{n-2}} {2}-\Delta_{n-1} -\frac {\Delta_{n-2}} {10}}{\Delta_{n-1}}+1 \\
    &\geq&\frac{1}{4} \frac{\Delta_{n-2}}{ \Delta_{n-1}}> s,\label{G318}
 \end{eqnarray}
 where $|T_\kappa|=b-a+1$ for $T_\kappa=[a,b]$.
For such $ T_\kappa$,
 a similar estimate to (\ref{G317}) gives
\begin{eqnarray}
\nonumber
  \sum_{j\in T_\kappa}\ln
|\sin  \pi\hat\theta_j|  &>&   |T_\kappa|\ln \frac {|T_\kappa|}{q_n} -C s  \\
 &>&    |T_\kappa|\ln \frac{s}{q_n}-Cs. \label{G319}
\end{eqnarray}
 If $T_\kappa$ does not   contain any $j$ with $ \|  \hat
\theta_j\|_{\mathbb{R}/\mathbb{Z}}<\frac {\Delta_{n-2}} {10} $, then by (\ref{G33})
\begin{eqnarray}
\nonumber
  \sum_{j\in T_\kappa}\ln
|\sin  \pi\hat\theta_j|  &>&  -|T_\kappa|\ln q_{n-1}-C|T_\kappa| \\
    &>&    |T_\kappa|\ln \frac{s}{q_n}-C|T_\kappa|,\label{G320}
\end{eqnarray}
 since $ s<\frac{q_{n}}{q_{n-1}}$.
 \par
By (\ref{G319}) and (\ref{G320}), one has
\begin{equation}\label{G321}
  \sum_{j\in J_1}\ln
|\sin  \pi\hat\theta_j|\geq  2s\ln \frac{s}{q_n} -C  s .
\end{equation}
Similarly,
\begin{equation}\label{G322}
  \sum_{j\in J_2}\ln
|\sin  \pi\hat\theta_j|\geq  4s\ln \frac{s}{q_n} -C  s .
\end{equation}
Putting  (\ref{G312}), (\ref{G321}) and (\ref{G322}) together, we have
\begin{equation}\label{G323}
\Sigma_+ > -6sq_{n-1}\ln 2+ 6s\ln \frac{s}{q_n} -Cs\ln q_{n-1}.
\end{equation}
We are now in the position to estimate $\Sigma_-$. Following the  discussion of $\Sigma_+$ , we have the similar estimate,
\begin{equation}\label{G324}
\Sigma_- > -6sq_{n-1}\ln 2+ 6s\ln \frac{s}{q_n} -Cs\ln q_{n-1}.
\end{equation}
 In order to avoid repetition, we omit the proof of (\ref{G324}).
\par
By  (\ref{G39}), (\ref{G323}) and (\ref{G324}), one obtains
$$ \sum _{j \in I_1  \cup I_2,j\neq i}\ln|\cos2\pi \theta_i-\cos2\pi \theta_j| \;\;\;\;\;\;\;\;\;\;\;\;\;\;\;\;\;\;\;\;\;\;\;\;\;\;\;\;\;\;\;\;\;\;\;\;\;\;\;\;\;\;\;\;\;\;\;\;\;\;\;\;\;\;\;\;$$
\begin{equation}\label{G325}
  > -6sq_{n-1}\ln 2+ 12s\ln \frac{s}{q_n} -Cs\ln q_{n-1}.
\end{equation}
Combining with (\ref{G38}),
 we have for any $\varepsilon>0$,
\begin{equation}\label{G326}
        \max_{ i\in I_1\cup I_2} \prod_{j \in I_1\cup I_2,j \neq i } \frac{|\cos2\pi a-\cos2\pi\theta_{j }|}
        {|\cos2\pi\theta_i-\cos2\pi\theta_{j }|}<e^{( 6sq_{n-1}-1)(-2\ln (s /q_n)/q_{n-1}+\varepsilon  )  },
      \end{equation}
      for $n $ large enough.$\qed$
\begin{theorem}\label{Avila1}
  Suppose $y$ is non-resonant. Let $s $   be
the largest positive integer  such that $4sq_{n-1}-p+1\leq \text{dist}(y,\{ \ell q_n\}_{\ell\geq0})\equiv y_0$. Then for
any $\varepsilon>0$ and   sufficiently large $n$,
$ y$ is $  (\ln \lambda+18\ln (s  q_{n-1}/q_n)/q_{n-1}-\varepsilon,6sq_{n-1}-1)$-regular if $\ln \lambda>2\beta $. In particular,
$ y$ is $  (\ln \lambda-2\beta-\varepsilon,6sq_{n-1}-1)$-regular.
\end{theorem}
\textbf{Proof}: Theorem \ref{Avila1}   can be derived from Lemma \ref{Le32} directly.
See the proof of Lemma 9.4 in \cite{AJ1} (p.24) for details.
\section{ Regularity for  resonant  $y$}
In this section, we mainly concern the regularity  for   resonant $y$.
If $b_n\leq y<b_{n+1}$ is  resonant, by the definition of resonance, there exists some  positive integer $ \ell$ with   $1 \leq \ell\leq q_{n+1}^{8/9}/q_n$
such that   $|y-\ell q_n|\leq b_n$.  Fix the positive integer $\ell$   and  let $s$ be the
largest positive integer such that $7sq_{n-1}\leq  q_n+p-1 $. Set $ I_1, I_2\subset \mathbb{Z}$ as follows
\begin{eqnarray*}
  I_1 &=& [-4sq_{n-1} ,-1] ,\\
  I_2 &=& [\ell q_n - 3sq_{n-1}  ,  \ell   q_n + 3sq_{n-1}  -1],
\end{eqnarray*}
and let $\theta_j=\theta+j\alpha$ for $j\in I_1\cup I_2$, the set $\{\theta_j\}_{j\in I_1\cup I_2}$
consists of $ 10 sq_{n-1}$ elements.
\par
  We will use the following three steps  to establish the regularity for $y$.
      \textbf{ Step 1}:  we  set up the $ \frac{7 }{5}\beta+\varepsilon$-uniformity of $\{\theta_j\}$  for any $\varepsilon>0$. By Lemma $ \ref{Le22}$, there exists some  $j_0\in I_1\cup I_2$
    such that
      $ \theta_{j _0}\notin  A_{10 sq_{n-1}-1,\ln\lambda- \frac{7 }{5}\beta - \varepsilon }$ for any  $\varepsilon>0$.   \textbf{ Step 2}: we   show that
      $\forall j \in I_1, \theta_j \in   A_{10 sq_{n-1}-1, \ln\lambda- \frac{7 }{5}\beta- \varepsilon }$ if $\lambda>e^{7\beta}$.  Thus there exists
       $ \theta_{j_0} \notin  A_{10 sq_{n-1}-1,\ln\lambda- \frac{7 }{5}\beta - \varepsilon}$ for some $j_0 \in I_2$.
        \textbf{ Step 3}:
       we establish the regularity for $y$.
      \par
      We start with the  \textbf{Step 1}.
       \begin{lemma}  \label{Le41}
For any     $ \varepsilon>0$,   the set $\{\theta_j\}_{j\in I_1\cup I_2}$
is $(\frac{7 }{5}\beta+\varepsilon)$-uniform  if  $n$ is sufficiently large.
\end{lemma}
\textbf{Proof:}
Notice that for any $ i\in I_1\cup I_2$,  there is at most one  $  \tilde{i}\in I_1\cup I_2$
such that $|i-\tilde{i}|=\ell q_n$.
It is easy to check
\begin{equation}\label{G41}
   \ln|\sin\pi(i-\tilde{i})\alpha| = \ln|\sin( \pi \ell q_n\alpha)| >-\ln q_{n+1}-C,
\end{equation}
since    $\Delta_n\geq \frac{1}{2q_{n+1}} $.
 If $j\neq i,\tilde{i}$ and $j\in I_1\cup I_2$, then $j-i=r+m_jq_n$ with $1\leq| r|<q_n$ and $|m _j|\leq \ell +2$. Thus
 by $(\ref{G32})$ and $(\ref{G33})$,
 $$||r\alpha||_{\mathbb{R}/ \mathbb{Z}}\geq \Delta_{n-1}  $$
 and
 \begin{eqnarray}
 \nonumber
    \min_{ {j\in I_1\cup I_2} {j\not=i,\tilde{i}}} || (j-i)\alpha||_{\mathbb{R}/ \mathbb{Z}} & > &    ||r\alpha||_{\mathbb{R}/ \mathbb{Z}}- (\ell +2) \Delta_n   \\
   & >& \frac{\Delta_{n-1} }{2}, \label{G42}
 \end{eqnarray}
since $ (\ell +2) \Delta_n  <\frac{\Delta_{n-1} }{2}$ for $n$ large enough.
\par
Similarly,
for any $ i\in I_1\cup I_2$,  there is at most one  $   \bar{i} \in I_1\cup I_2$
such that $|i+\bar{i} +p|=\ell q_n$. We also have
\begin{equation}\label{XXG43}
     \ln|\sin\pi(2\theta+(i+\bar{i})\alpha) |   >-\ln q_{n+1}-C,
\end{equation}
and
\begin{equation}\label{XXG44}
  \min_{ {j\in I_1\cup I_2} {j\not=\bar{i},  }} ||  2\theta+(j+i)\alpha ||_{\mathbb{R}/ \mathbb{Z}}   > \frac{\Delta_{n-1} }{2}.
\end{equation}
Replacing  (\ref{G316}) with (\ref{G41}) and  (\ref{G42})  and following the proof of Lemma \ref{Le32},
one has
\begin{equation}\label{XXXG45}
\Sigma_- > -10sq_{n-1}\ln 2+ 10s\ln \frac{s}{q_n} -\ln q_{n+1}-Cs\ln q_{n-1},
\end{equation}
since  there exists at most one  term satisfies (\ref{G41}).
\par
Similarly, Replacing  (\ref{G316}) with   (\ref{XXG43}) and (\ref{XXG44}), one has
\begin{equation}\label{XXXG46}
\Sigma_+> -10sq_{n-1}\ln 2+ 10s\ln \frac{s}{q_n} -\ln q_{n+1}-Cs\ln q_{n-1}.
\end{equation}
By (\ref{G38}), (\ref{XXXG45}) and (\ref{XXXG46}),
we have for any $\varepsilon_0>0$,
\begin{equation}\label{G43}
        \max_{ i\in I_1\cup I_2} \prod_{j \in I_1\cup I_2,j \neq i } \frac{|x-\cos2\pi\theta_{j }|}
        {|\cos2\pi\theta_i-\cos2\pi\theta_{j }|}<e^{ 10sq_{n-1}(-2\ln (s /q_n)/q_{n-1}+ 2\frac{\ln q_{n+1}}{10sq_{n-1}}+\varepsilon_0  )  },
      \end{equation}
      if  $n $ is large enough.
      \par
      By the definition of $s$ and noting that $\beta=  \limsup_{n\rightarrow\infty}\frac{\ln q_{n+1}}{q_n}$,  one has
      \begin{equation}\label{G44}
       -2\ln (s /q_n)/q_{n-1}+  \frac{\ln q_{n+1}}{5sq_{n-1}}< \frac{7}{5} \beta+\varepsilon_0,
      \end{equation}
        for $n $  large enough.
        Combining (\ref{G43}) with (\ref{G44}), we obtain
        \begin{equation}\label{G45}
        \max_{ i\in I_1\cup I_2} \prod_{j \in I_1\cup I_2,j \neq i } \frac{|x-\cos2\pi\theta_{j }|}
        {|\cos2\pi\theta_i-\cos2\pi\theta_{j }|}<e^{ (10sq_{n-1}-1)( \frac{7}{5} \beta+3\varepsilon_0  )  }.
      \end{equation}
       By the arbitrariness of $\varepsilon_0 $,   we complete the proof.$\qed$
       \par
      Now, we are in the position to undertake \textbf{ Step 2}.
\begin{lemma}\label{Le42}
 Assume $ y\in [- 2q_n, 2q_n ]$ and let $d= \text{dist}(y,\{j q_n\}_{j\geq0})>\frac{1}{100}q_n$.
 Then, for any  $\varepsilon>0$,
\begin{equation}\label{G46}
    |\phi(y)|<\exp(-(\ln \lambda-\varepsilon)d)
\end{equation}
 if $n$ is sufficiently large.
\end{lemma}
\textbf{Proof:} Using Theorem \ref{Avila1} and  block-resolvent expansion,  it is easy to obtain   Lemma \ref{Le42}.
   See the proof of     Lemma 3.4 in $ \cite{LIU}$ for details.
\begin{theorem}\label{Th43}
For any $  \varepsilon>0$ and   any $b\in[- 9s q_{n-1},- 5sq_{n-1}]$,
we have $\theta+ ( b+ 5sq_{n-1}  -1)\alpha\in A_{ 10sq_{n-1}-1,\frac{4}{5} \ln \lambda+\varepsilon}$ if $n$ is large enough. That is  for all $j\in I_1$,
$\theta_j\in A_{ 10sq_{n-1}-1,\frac{4}{5} \ln \lambda +\varepsilon}$.
\end{theorem}
\textbf{Proof:} For any    $b\in[- 9s q_{n-1},- 5sq_{n-1}]$, let $b_1=b-1$ and $b_2=b+ 10sq_{n-1}-1$.
Applying  Lemma $\ref{Le42}$,  one has, for any $ \varepsilon_0$,
 \begin{equation}\label{G411}
 |\phi(b_1)|<e^{-(\ln \lambda-\varepsilon_0)|q_{n}+b|} ,|q_{n}+b|>\frac{q_n}{100},
 \end{equation}
 and
  \begin{equation}\label{G412}
  |\phi(b_2)|\leq \left\{
                \begin{array}{ll}
                  e^{-(\ln\lambda-\varepsilon_0)(10sq_{n-1}+ b)}, &  {   b \in [-9sq_{n-1}, \frac{q_n}{2}-10sq_{n-1} ] ;} \\
                  e^{ -(\ln\lambda-\varepsilon_0) (q_n-10sq_{n-1}- b) }, &  {    b \in [  \frac{q_n}{2}-10sq_{n-1} ,-5sq_{n-1}].}
                \end{array}
              \right.
  \end{equation}
  By the definition of $s$, (\ref{G411}) and (\ref{G412})
  become
  \begin{equation}\label{G413}
 |\phi(b_1)|<e^{-(\ln \lambda-2\varepsilon_0)| 7sq_{n-1}+b|} ,|q_{n}+b|>\frac{q_n}{100},
 \end{equation}
 and
  \begin{equation}\label{G414}
  |\phi(b_2)|\leq \left\{
                \begin{array}{ll}
                  e^{-(\ln\lambda-\varepsilon_0)(10sq_{n-1}+ b)}, &  {   b \in [-9sq_{n-1}, \frac{q_n}{2}-10sq_{n-1} ] ;} \\
                  e^{ -(\ln\lambda-\varepsilon_0) ( -3sq_{n-1}- b) }, &  {    b \in [  \frac{q_n}{2}-10sq_{n-1} ,-5sq_{n-1}].}
                \end{array}
              \right.
  \end{equation}
  \par
In ($ \ref{G29}$), let $x=0$  and  $I=[b,b+ 10sq_{n-1}-2]$, we get for $n$ large enough,
\begin{equation}\label{G415}
 |G_I(0,b)|> \left\{
                \begin{array}{ll}
                 e^{ (\ln\lambda-3\varepsilon_0)|7sq_{n-1}+ b|}, &  { |q_{n}+b|>\frac{q_n}{100} ;} \\
                  e^{     -  \varepsilon_0sq_{n-1}   }, &  {   |q_{n}+b|\leq \frac{q_n}{100},}
                \end{array}
              \right.
  \end{equation}

or
\begin{equation}\label{G416}
 |G_I(0,b+ 10sq_{n-1}-2)|> \left\{
                \begin{array}{ll}
                  e^{ (\ln\lambda-3\varepsilon_0)(10sq_{n-1}+ b)}, &  {   b \in [-9sq_{n-1}, \frac{q_n}{2}-10sq_{n-1} ] ;} \\
                  e^{ (\ln\lambda-3\varepsilon_0) ( -3sq_{n-1}- b) }, &  {    b \in [  \frac{q_n}{2}-10sq_{n-1} ,-5sq_{n-1}],}
                \end{array}
              \right.
  \end{equation}
  since  $\phi(0)=1$ and $|\phi(k)|\leq (1+|k|)^C $.
\par
By ($\ref{G25}$), ($\ref{G26} $) and ($ \ref{G27}$),
\begin{flushleft}
   $|  Q_{ 10sq_{n-1}-1}(\cos2\pi (\theta+(b+\frac{10sq_{n-1}-2}{2})\alpha)|$
\end{flushleft}
$$ =|P_{ 10sq_{n-1}-1}(\theta+b\alpha)|\;\;\;\;\;\;\;\;\;\;\;\;\;\;\;\;\;\;\;\;\;\;\;\;\;\;\;\;\;\;\;\;\;\;\;\;\;\;\;\;\;\;\;\;\;\;\;\;\;\;\;\;
\;\;\;\;\;\;\;\;\;\;\;\;\;\;\;\;\;\;\;\;\;\;$$
$$<\min\{|G_I(0,b)|^{-1}e^{(\ln \lambda+\varepsilon_0)(b+10sq_{n-1}-2)}, |G_I(0,b+ 10sq_{n-1}-2)|^{-1}e^{-(\ln \lambda+\varepsilon_0)b}\}$$
$$ <e^{ (\frac{4}{5}\ln \lambda+4\varepsilon_0)10sq_{n-1}}. \;\;\;\;\;\;\;\; \;\;\;\;\;\;\;\;\;\;\;\;\;\;\;\;\;\;\;\;\;\;\;\;\;\;\;\;\;\;\;\;\;\;\;\;\;\;\;\;\;\;\;\;\;\;\;\;\;
\;\;\;\;\;\;\;\;\;\;\;\;\;\;\;\;\;\;\;\;$$
This implies  $\theta+ ( b+ 5sq_{n-1}  -1)\alpha\in A_{ 10sq_{n-1}-1,\frac{4}{5} \ln \lambda+ 4\varepsilon_0}$.
By the arbitrariness of $\varepsilon_0 $,
we have $\theta+ ( b+ 5sq_{n-1}  -1)\alpha\in A_{ 10sq_{n-1}-1,\frac{4}{5} \ln \lambda+\varepsilon}$ for
any $b\in[- 9s q_{n-1},- 5sq_{n-1}]$.  $\qed$
\par
Finally, we will finish the \textbf{Step 3}.
\begin{theorem}\label{Th44}
 For any
$ \varepsilon>0$  such that $t= (\ln\lambda-7\beta-\varepsilon)>0$,
  $y$ is $(t,10sq_{n-1}-1 )$-regular if $n$ is  large enough  (or  equivalently  $y$ is large enough).
\end{theorem}
\textbf{Proof:} Let $y>0$ be resonant. By hypothesis  $ y=\ell q_n+r$,
with $ 0\leq |r|\leq q_n^{\frac{8}{9}}$ and $1\leq \ell \leq q_{n+1}^{\frac{8}{9}}/q_n.$
\par
 By Lemma $\ref{Le22}$ and Lemma  $\ref{Le41}$ (let $\varepsilon=\varepsilon_0/2$ in Lemma  $\ref{Le41}$), for any $\varepsilon_0>0$, there exists some   $ j\in I_1\cup I_2  $ such that
 $\theta_j\notin A_{10sq_{n-1}-1,\ln \lambda-\frac{7}{5}  \beta -\varepsilon_0}$.
 By Theorem  $ \ref{Th43}$ and noting that  $\ln \lambda>  7\beta$ (i.e.\;$ \frac{4}{5}\ln \lambda<(\ln \lambda-\frac{7}{5}  \beta )$),
we have $ \theta_j\in A_{10sq_{n-1}-1,\ln \lambda- \frac{7}{5} \beta -\varepsilon_0} $
for all $j\in I_1$ and sufficiently small $\varepsilon_0$. Thus, there exists some $j_0\in I_2$ such that $\theta_{j_0}\notin A_{10sq_{n-1}-1,\ln \lambda- \frac{7}{5} \beta -\varepsilon_0}$.
Set $I=[j_0-5sq_{n-1} +1,j_0 +5sq_{n-1} -1]=[x_1,x_2]$.
By ($\ref{G25}$), ($\ref{G26} $) and ($ \ref{G27}$) again,  we have
$$|G_I(y,x_i)|<e^{(\ln \lambda+\varepsilon_0)(10sq_{n-1}-2-|y-x_i|)-  10sq_{n-1} (\ln\lambda-\frac{7}{5}  \beta -\varepsilon_0)}.$$
By a simple computation
$$|y-x_i|\geq (2sq_{n-1}- q_n^{\frac{8}{9}})>(\frac{1}{5}-\varepsilon_0) 10sq_{n-1}, $$
 therefore,
 $$|G_I(y,x_i)|<e^{-|y-x_i|( \ln \lambda-7 \beta -\varepsilon )},$$
 where
$\varepsilon  = 20\varepsilon_0$. Let $t=
\ln \lambda- 7\beta -\varepsilon>0 $, then for $n$ large
enough, $y$ is $(t,10sq_{n-1}-1 )$-regular.$\qed$
\section{Proof of Theorem  \ref{Main Theorem}}
Now that  the regularity for   $y$ is  established, we will  use  block resolvent expansion  to prove Theorem $\ref{Main Theorem}$.
\par
\textbf{Proof of     Theorem \ref{Main Theorem}.}

  Give some $k$ with  $   k>q_n$ and $n$ large enough.
$\forall y \in [ q_n^{\frac{8}{9}}, 2k]$,   let $\varepsilon=\varepsilon_0$  in   Theorem $ \ref{Avila1}$  and $ \ref{Th44}$,  then there exists an interval $ I(y)=[x_1,x_2]\subset
[ -4k,4k]$ with
 $y\in I(y)$ such that
\begin{eqnarray}\label{G51}
\nonumber
    \text{dist}(y,\partial I(y)) &>& \frac{ 1}{7} |I(y)| \geq\min{\{\frac{6sq_{n-1}-1}{7},\frac{10sq_{n-1}-1}{7} \}}  \\
    &\geq&  \frac{1}{2}q_{n-1}
\end{eqnarray}
and
\begin{equation}\label{G52}
  |G_{I(y)}(y,x_i)| < e^{-(\ln\lambda-7\beta-\varepsilon_0)|y-x_i|},\;i=1,2.
\end{equation}
 Denote by $ \partial I(y)$    the boundary of the interval $I(y)$. For $z  \in  \partial I(y)$,
let $z' $ be the neighbor of $z$, (i.e., $|z-z'|=1$) not belonging to $I(y)$.
\par
If $x_2+1< 2k$ or  $x_1-1>  b_n=q_n^{\frac{8}{9}}$,
we can expand $\phi(x_2+1)$ or $\phi(x_1-1)$ as ($ \ref{G29}$). We can continue this process until we arrive to $z$
such that $z+1\geq 2k$ or  $z-1\leq  b_n$, or the iterating number reaches
$[\frac{2k}{q_{n-1}}]$.
\par  By (\ref{G29}),
\begin{equation}\label{G53}
   \phi(k)=\displaystyle\sum_{s ; z_{i+1}\in\partial I(z_i^\prime)}
G_{I(k)}(k,z_1) G_{I(z_1^\prime)}
(z_1^\prime,z_2)\cdots G_{I(z_s^\prime)}
(z_s^\prime,z_{s+1})\phi(z_{s+1}^\prime),
\end{equation}
where in each term of the summation we have
$b_n+1<z_i<2k-1$, $i=1,\cdots,s,$ and
  either $z_{s+1} \notin [b_n+2, 2k-2]$, $s+1 < [\frac{2 k}{q_{n-1 }}]$; or
$s+1= [\frac{2k}{q_{n-1}}]$.
\par
 If $z_{s+1} \notin [b_n+2, 2k-2]$, $s+1 < [\frac{2k}{q_{n-1 }}]$, by  ($\ref{G52}$) and noting that $|\phi(z_{s+1}^\prime)|\leq(1+|z_{s+1}^\prime|)^C\leq k^C $, one has
\begin{eqnarray}\label{G54}
\nonumber
 & | G_{I(k)}(k,z_1) G_{I(z_1^\prime)}
(z_1^\prime,z_2)\cdots G_{I(z_s^\prime)}
(z_s^\prime,z_{s+1})\phi(z_{s+1}^\prime)|\;\;\;\;\;\;\;\;\;\;\;\;\;\;\;\;\;\;\;\;\;\;\;\;\;\;\;\;\;\;\;\;\;\;\\
\nonumber
&\le e^{-(\ln\lambda- 7\beta-\varepsilon_0)(|k-z_1|+\sum_{i=1}^{s}|z_i^\prime-z_{i+1}|)}
k ^C\;\;\;\;\;\;\;\;\;\;\;\;\;\;\;\;\;\;\;\;\;\;\;\;\;\;\;\;\;\;\;\;\;\;\; \;\;\;\; \; \\
\nonumber
&\le e^{-(\ln\lambda-  7\beta-\varepsilon_0 )(|k-z_{s+1}|-(s+1))}k ^C\;\;\;\;\;\;\;\;\;\;\;\;\;\;\;\;\;\;\;\;\;\;\;\;\;\;\;\;\;\;\;\;\;\;\;\;\;\;\;\;\;\;\;\;\;\; \\
   &\le \max\{e^{-(\ln\lambda-7 \beta-\varepsilon_0 )(k- b_n -4-\frac{ 2k}{q_{n-1}})}k^C ,
e^{-(\ln\lambda-  7\beta-\varepsilon_0)(2k-k -4-\frac{ 2k}{q_{n-1}})}k^C \}.
\end{eqnarray}
If $s+1= [\frac{2k}{q_{n-1}}] ,$
using  ($\ref{G51}$) and ($\ref{G52}$),     we obtain
\begin{equation}\label{G55}
     | G_{I(k)}(k,z_1) G_{I(z_1^\prime)}
(z_1^\prime,z_2)\cdots G_{I(z_s^\prime)}
(z_s^\prime,z_{s+1})\phi(z_{s+1}^\prime)|\le  e^{-(\ln\lambda-7\beta-\varepsilon_0) {\frac{q_{n-1}}{2}} [\frac{2 k}{q_{n-1}}]} k^C.
\end{equation}

Finally,   notice that the total number of terms in ($ \ref{G53}$)
is  at most  $2^{[\frac{2k}{q_{n-1}}]}$. Combining with
($\ref{G54}$) and
($\ref{G55}$),  we obtain
\begin{equation}\label{G56}
|\phi(k)|\le   e^{-(\ln\lambda-7\beta-2\varepsilon_0-\varepsilon_0\ln\lambda) k }
\end{equation}
for large enough $n$ (or equivalently  large enough $k$ ).  By the arbitrariness of $\varepsilon_0 $, we have for any $\varepsilon>0$,
\begin{equation}\label{G57}
|\phi(k)|\le   e^{-(\ln\lambda-7\beta-\varepsilon) k }  \text{  for } k \text{  large enough}.
\end{equation}
\par
For $k<0$, the proof is similar.
Thus  for any $\varepsilon>0$,
\begin{equation}\label{G58}
|\phi(k)|\le   e^{-(\ln\lambda-7\beta-\varepsilon) |k| }  \text{  if } |k| \text{  is large enough}.
\end{equation}
We finish the proof of    Theorem $\ref{Main Theorem}$.

\par
           \begin{center}
           
             \end{center}
  \end{document}